\documentclass[reqno]{amsart}
\usepackage{amsmath}
\usepackage{amsthm}
\usepackage{amssymb}
\usepackage{color,soul,comment}

\newcommand{\ds}{\displaystyle}
\newcommand{\reals}{\mathbb{R}}
\newcommand{\integers}{\mathbb{Z}}
\newcommand{\DN}{\mathsf{DN}}
\newcommand{\diag}{\text{diag}}

\theoremstyle{plain}
\newtheorem{theorem}{Theorem}[section]
\newtheorem{lemma}{Lemma}[section]

\theoremstyle{definition}

\title[The Critical Exponent for Doubly Nonnegative Matrices]{The Critical Exponent for Continuous Conventional Powers of Doubly Nonnegative Matrices}
\author{Charles R. Johnson}
\address{Charles R. Johnson, The College of William \& Mary}
\author{Brian Lins}
\address{Brian Lins, Hampden-Sydney College \footnote{Corresponding author}}
\author{Olivia Walch}
\address{Olivia Walch, The College of William \& Mary}

\thanks{This work was partially supported by NSF grant DMS-0751964} 

\begin{document}

\begin{abstract}
We prove that there exists an exponent beyond which all continuous conventional powers of $n$-by-$n$ doubly nonnegative matrices are doubly nonnegative.  We show that this critical exponent cannot be less than $n-2$ and we conjecture that it is always $n-2$ (as it is with Hadamard powering).  We prove this conjecture when $n<6$ and in certain other special cases.  We establish a quadratic bound for the critical exponent in general.  
\end{abstract}
\maketitle
\section{Introduction}
An $n$-by-$n$ real symmetric matrix is called doubly nonnegative ($\DN$) if it is both positive semi-definite and entry-wise nonnegative.  Continuous powers (with exponent at least 0) of a $\DN$ matrix $A=(a_{ij})$ are well defined under both conventional multiplication, $A^t$ (defined in the primary way \cite{HornJohnson2} via the spectral decomposition), and under Hadamard multiplication, $A^{(t)} = (a^t_{ij})$.  In both cases, positive integral powers are well-known to remain $\DN$ (in the Hadamard case, because of Schur's observation \cite{HornJohnson1}).  But, in both cases, there is a natural question about powers between positive integral powers.  

In the Hadamard case, entry-wise nonnegativity remains clear, but the nonnegativity of the quadratic form is a question.  In the continuous case, the nonnegativity of the quadratic form of all powers is clear, but the entry-wise nonnegativity of continuous powers is not.  In both cases, this suggests the definition (whose appropriateness must be proven) of a \textbf{critical exponent}.  The Hadamard (conventional) critical exponent is the least value $m$, for a given $n$, such that $A^{(t)}$ ($A^t$) is $\DN$ for all $t \geq m$ and all $A$ that are $\DN$.  Of course, it is not immediately clear that either critical exponent is finite, but both are.  

The Hadamard critical exponent has been studied \cite{FH} in conjunction with the (then) Bieberbach conjecture and, interestingly, been shown to be $n-2$.  The natural guess of one is false, though it is correct for the $\DN$ matrices that are also inverse $M$-matrices \cite{Chen}.  

Here we take up the issue of the conventional $\DN$ critical exponent, which is of interest not only by analogy, but also because of the fundamental interest in entry-wise nonnegativity and dynamical systems.  We show that the conventional critical exponent is also finite, that it is at least $n-2$ (because of tridiagonal $\DN$ matrices), and we give low-coefficient quadratic upper bounds for it.  We conjecture, interestingly, that the conventional critical exponent is also $n-2$, though there appears to be little technical relation between the two.  This conjecture is proven for $n<6$.  Facts about exponential polynomials are exploited and some new techniques, based upon a combinatorial matrix built from eigenvectors, are developed to do this.  Additional observations, involving entry-wise matrix inequalities for powers, etc., are also made.   

\section{Background}

Any symmetric matrix $A \in M_n(\reals)$ can be decomposed as $A = \lambda_1 x_1 x_1^T + ... + \lambda_n x_n x_n^T$ where the set $\{x_1, \ldots ,x_n\}$ is an orthonormal basis of eigenvectors and for each $x_i$, $\lambda_i$ is the corresponding eigenvalue.  If $A$ is positive semi-definite, then for $t \in \reals \backslash \{ 0 \}$, $A^t$ is defined by
$$A^t = \lambda_1^t x_1 x_1^T + \cdots +\lambda_n^t x_n x_n^T.$$
Each entry of $A^t$ has the form 
$$(A^t)_{ij} = \lambda_1^t (x_1 x_1^T)_{ij} + \cdots \lambda_n^t (x_n x_n^T)_{ij}.$$
Any function of the form 
$$\varphi(t) = \alpha_1 e^{\beta_1 t} + \cdots \alpha_n e^{\beta_n t}$$ 
is an \textbf{exponential polynomial}.  In particular, if $A$ is positive semi-definite, then each entry of $A^t$ is an exponential polynomial in $t$.  The following version of Descartes' rule for exponential polynomials is well known and appears as an exercise in \cite{PoSz}.  

\begin{lemma} \label{descartes} 
Let $\varphi(t) =\sum_{i = 1}^k \alpha_i e^{\beta_i t}$ be a real exponential polynomial such that each $\alpha_i \neq 0$ and $\beta_1 > \beta_2 > \ldots > \beta_n$.  The number of real roots of $\varphi(t)$, counting multiplicity, cannot exceed the number of sign changes in the sequence of coefficients $\{\alpha_1, \alpha_2, \ldots, \alpha_n\}$.
\end{lemma}

Lemma \ref{descartes} leads immediately to the following theorem, which proves the existence of a critical exponent for continuous conventional powers of doubly nonnegative matrices.  

\begin{theorem} \label{existence}
There is a function $m(n)$ such that for any $n$-by-$n$ doubly nonnegative matrix $A$, $A^t$ is doubly nonnegative for $t \geq m(n)$.  
\end{theorem}

\begin{proof}
Let $A$ be an $n$-by-$n$ doubly nonnegative matrix.  Since $A$ is nonnegative, so is $A^k$ for all positive integers $k$.  If $A$ is nonnegative for all $t \in [m,m+1]$, where $m \in \integers$, then it follows from repeated multiplication by $A$ that $A^t$ is nonnegative for all $t \geq m$.  Suppose that $A^t$ has a negative entry for some $t \in [m,m+1]$, then the exponential polynomial corresponding to that entry must have at least two roots in the interval $[m,m+1]$.  By Lemma \ref{descartes}, the maximum number of roots each entry may possess depends on $n$.  It follows that there is a constant $m(n)$ such that $A^t$ is nonnegative for all $t > m(n)$.  
\end{proof}

Let $A$ be any $n$-by-$n$ doubly nonnegative matrix.  Corresponding to the matrix $A$, we define a matrix $W = [w_{ij}]$ where $w_{ij}$ equals the number of sign changes in the sequence of coefficients of the exponential polynomial $(A^t)_{ij}$ arranged in decreasing order of the corresponding eigenvalues.  We refer to any matrix defined this way as the \textbf{sign change matrix} for $A$.  By Lemma \ref{descartes}, each entry $w_{ij}$ of a sign change matrix gives an upper bound on the number of real zeros of the corresponding exponential polynomial $(A^t)_{ij}$, counting multiplicity.  The following lemma gives some restrictions on the structure of a sign change matrix.  

\begin{lemma} \label{Wmatrix}
Let $A$ be an $n$-by-$n$ doubly nonnegative matrix.  If $W$ is the sign change matrix corresponding to $A$, then every diagonal entry of $W$ is zero, every row and column of $W$ contains at most one entry equal to $n-1$, and the remaining entries of $W$ are at most $n-2$.  
\end{lemma}
\begin{proof}
Since $A$ is symmetric, there is an orthogonal matrix $U$ such that $A = U D U^T$ where $D$ is the diagonal matrix $D = \diag(\lambda_1, \lambda_2, ..., \lambda_n)$ and $U = [u_{ij}]$.  The $i,j$-entry of $A^t$ is given by
$$(A^t)_{ij} = e_i^T U D^t U^T e_j = (U^T e_i)^T D^t (U^T e_j) =$$
$$= u_{i1} u_{j1} \lambda_1^t + \cdots + u_{in} u_{jn} \lambda_n^t.$$
Note that the coefficients of the the exponential polynomial $(A^t)_{ij}$ are given by the Hadamard product of the $i^{\text{th}}$ row of $U$ with the $j^{\text{th}}$ row of $U$.  Since $U$ is an orthogonal matrix, no two rows of $U$ can have the same sign pattern.  Therefore, in any given column or row of $A^t$, only one entry can have an exponential polynomial with $n-1$ sign changes.  The remaining entries in the columns may have at most $n-2$ sign changes.  Furthermore, the diagonal entries of $A^t$ have exponential polynomials with all nonnegative coefficients, so there are zero sign changes.  
\end{proof}

\begin{lemma} \label{components}
Let $A$ be an invertible doubly nonnegative matrix with sign change matrix $W=[w_{ij}]$.  Let $T^-_{ij} = \{t > 1 : (A^t)_{ij} < 0 \}$.  Then the maximum number of connected components of $T^-_{ij}$ is 
$$\begin{cases}
\lfloor (w_{ij} - 1)/2 \rfloor & \text{if}~w_{ij} > 0 \\
0 & \text{if}~w_{ij}=0.
\end{cases}$$
\end{lemma}
\begin{proof}
By Lemma \ref{descartes}, the maximum number of real roots of the exponential polynomial $(A^t)_{ij}$ is given by $w_{ij}$.  Since $A$ is invertible, the exponential polynomials defining the entries of $A^t$ when $t>0$ still agree with $A^t$ at $t=0$.  Since $A^0$ is the identity matrix, the exponential polynomial $(A^t)_{ij}$ has at most $w_{ij} -1$ roots in the interval $[1,\infty)$ when $i \neq j$.  

Each of the connected components of $T_{ij}^-$ is bounded because $A^k$ is nonnegative for all positive integers $k$.  The endpoints of these components are roots of the exponential polynomial $(A^t)_{ij}$.  If two adjacent connected components of $T^-_{ij}$ share an endpoint, that endpoint must be a root of degree at least two.  Counting multiplicity, the number of real roots of $(A^t)_{ij}$ with $t\geq 1$ must therefore be at least double the number of connected components of $T^-_{ij}$.

If $w_{ij}$ is zero, then the exponential polynomial $(A^t)_{ij}$ has all positive coefficients, so $T^-_{ij}$ is empty.  Note that $w_{ij} = 0$ whenever $i = j$ by Lemma \ref{Wmatrix}.  
\end{proof}

\begin{lemma} \label{primitive}
If $A \in M_n$ is irreducible and doubly nonnegative, then $A$ is primitive with index of primitivity at most $n-1$.  
\end{lemma}

\begin{proof}
Since $A$ is doubly nonnegative, every entry on the main diagonal is positive.  For an irreducible matrix $A$ with positive main diagonal, it is a routine exercise to verify that $A^{n-1}$ is entry-wise positive (see e.g., Lemma 8.5.5 in \cite{HornJohnson1}).  
\end{proof}

\section{Bounds for the Critical Exponent}

Using the tools developed in the previous section, we are now able to give an upper bound for the critical exponent.  

\begin{theorem} \label{crude}
The critical exponent $m(n)$ satisfies 
$$ m(n) \leq \begin{cases} \ds \frac{n^2-4n+5}{2} & \text{if} ~ n ~\text{odd} \\ \ds \frac{n^2-5n+8}{2} & \text{if} ~ n ~ \text{even}. \end{cases}$$
\end{theorem}

\begin{proof}
Suppose that $A$ is an $n$-by-$n$ doubly nonnegative matrix.  Assume for now that $A$ is irreducible and invertible.  Let $W$ denote the sign change matrix corresponding to $A$.  For each pair $i,j$, let $T_{ij}^- = \{t > 1: (A^t)_{ij} < 0 \}$.  Fix any column $j$ and let $T_j^- = \bigcup_{1 \leq i \leq n} T_{ij}^-$.  Using Lemma \ref{Wmatrix} and Lemma \ref{components} we arrive at the following upper bounds for the number of connected components of $T_j^-$.  

If $n$ is odd, then one entry in column $j$ with $w_{ij} = n-1$ can have up to $(n-2)/{2}$ connected components in $T_{ij}^-$.  The diagonal entry corresponds to $T_{jj}^- = \varnothing$.  The remaining entries correspond to sets $T_{ij}^-$ with up to ${(n-4)}/{2}$ connected components.  Therefore $T_j^-$ can have at most 
$$\frac{(n-2)}{2}+\frac{(n-2)(n-4)}{2} = \frac{n^2 - 5n + 6}{2}$$ 
connected components.  

If $n$ is even, then any entry in column $j$ with $w_{ij} = n-1$ or $n-2$ can have up to $(n-3)/{2}$ connected components in $T_{ij}^-$.  The diagonal entry corresponds to $T_{jj}^- = \varnothing$.  Therefore $T_j^-$ can have no more than 
$$\frac{(n-1)(n-3)}{2} = \frac{n^2 - 4n + 3}{2}$$ 
connected components.  

Since $A$ is nonnegative, so is $A^k$ for all positive integers $k$.  Therefore, each connected component of $T_j^-$ is contained an open interval $(m,m+1)$ for some positive integer $m$.  Let $k(n)$ equal the estimate given above for the number of connected components of $T^-_{j}$, that is
$$k(n) = \begin{cases} \ds \frac{n^2 - 4n + 3}{2}  & \text{if} ~ n ~\text{odd} \\ \ds \frac{n^2 - 5n + 6}{2}  & \text{if} ~ n ~\text{even}. \end{cases}$$
Suppose that $T_j^- \cap (m,m+1) = \varnothing$, for some $m \in \{1, \ldots, k(n) \}$.  This means that every entry in column $j$ of $A^t$ is nonnegative for all powers $t \in [m,m+1]$.  Using repeated left multiplication by $A$, we see that column $j$ of $A^t$ must be nonnegative for all $t \geq m$.  If $T_j^-$ has a connected component in each interval $(m,m+1)$, $m \in \{1, \ldots, k(n) \}$, then since $T^-_j$ has at most $k(n)$ connected components, it follows that $T_j^- \cap (k(n)+1,\infty) = \varnothing$.  Either way, the $j^{th}$ column of $A^t$ is nonnegative for all $t \geq k(n)+1$.  Since this applies to every column index $j \in \{1, \ldots, n\}$, we conclude that $A^t$ is doubly nonnegative for all $t \geq k(n)+1$.

Up until now, we have assumed that $A$ is both invertible and irreducible.  Suppose now that $A$ is reducible but still invertible.  Since $A$ is symmetric, it follows that there is a permutation matrix $P$ such that $PAP^T$ is a direct sum of smaller irreducible doubly nonnegative matrices and possibly 1-by-1 blocks containing zero.  Therefore the critical exponent of $A$ is bounded by the critical exponent of the these smaller blocks.  Since $k(n)$ is monotone, we see that $A^t$ is doubly nonnegative for all $t \geq k(n)+1$.  

Now suppose that $A$ is singular.  By continuity, $A^t$ cannot have a negative entry for any $t > k(n)+1$.  Therefore the critical exponent $m(n) \leq k(n)+1$.  
\end{proof}

The bound established in Theorem \ref{crude} is optimal for $n = 3, 4$ as we will show in Theorem \ref{lowerbound}.  We will show that it is not optimal when $n = 5$, by proving a sharper upper-bound in Section \ref{n=5}.

Using tridiagonal matrices, we can prove the following lower bound for the critical exponent of doubly nonnegative matrices.  

\begin{theorem} \label{lowerbound}
The critical exponent $m(n) \geq n-2$.  
\end{theorem}

\begin{proof}
Let $A$ be an invertible, irreducible, tridiagonal, doubly nonnegative matrix.  Since $A$ is irreducible, it must be primitive by Lemma \ref{primitive}.  Note that the $1,n$-entry of $A^t$ is zero for $t = 0, 1, 2, \ldots, n-2$.  By Lemma \ref{descartes}, the exponential polynomial $(A^t)_{1n}$ has at most $n-1$ zeros counting multiplicity.  Therefore $(A^t)_{1n} > 0$ for all $t > n-2$ and $(A^t)_{1n} < 0$ for all $t \in (n-3,n-2)$.  Thus $m(n) \geq n-2$.  
\end{proof}

\section{Critical Exponent when $n = 5$} \label{n=5}
In this section we prove that the critical exponent for 5-by-5 matrices doubly nonnegative matrices is $n-2$. 
\begin{theorem} \label{5by5}
The critical exponent $m(5) = 3$.  
\end{theorem}

For a doubly nonnegative matrix $A$, we define the \textbf{critical exponent of the $i,j$-entry of $A$} to be the least value of $m$ such that $(A^t)_{ij} \geq 0$ for all $t\geq m$.  

\begin{lemma} \label{Wbound12}
Let $A$ be an invertible doubly nonnegative matrix and let $W = [w_{ij}]$ be the sign change matrix corresponding to $A$.  The critical exponents of each entry in $A$ depend on the corresponding entry in $W$ as follows.  
\begin{enumerate}
\item If $w_{ij} = 0$ or $1$, then the critical exponent of the $i,j$ entry is 0.  
\item If $w_{ij} = 2$, then the critical exponent of the $i,j$ entry is 1. 
\end{enumerate}
\end{lemma}
\begin{proof}
Recall that the $i,j$-entry of the $W$ gives an upper bound on the number of real zeros of the exponential polynomial $(A^t)_{ij}$, counting multiplicity.  Note that $w_{ij} = 0$ if and only if $i = j$, and it is clear that the diagonal entries of $A^t$ are positive for all $t$ since $A$ is positive definite.  Since $A$ is invertible, $A^t$ is continuous at $t=0$ and $A^0 = I_n$.  Therefore, when $i \neq j$, the exponential polynomial $(A^t)_{ij}$ has a zero at $t=0$.  If $w_{ij}=1$, then we conclude that $(A^t)_{ij} > 0$ for all $t > 0$. Similarly, if $w_{ij} = 2$, then either the $(A^t)_{ij}$ is nonnegative after 0 (using only one of its two allotted zeros), or it is negative over an interval $(0,\epsilon)$ where $\epsilon < 1$ (using both zeros).  Thus the maximum exponent at which an entry corresponding to a 2 in the sign change matrix can be negative is $t=1$.  
\end{proof}

The only entries of $A^t$ that can be negative for $t>1$ correspond to the entries of $W$ that are larger than 2.  The following lemma addresses these entries.

\begin{lemma} \label{Wbound34}
Let $A$ be an invertible doubly nonnegative matrix and let $W=[w_{ij}]$ be the sign change matrix corresponding to $A$.  If a row (or column) of $W$ contains no entry greater than 4 and at most $M$ entries greater then 2, then the critical exponent of every entry in the corresponding row (or column) of $A$ is at most $M+1$.  
\end{lemma}  
\begin{proof}
Since $A$ is invertible, each exponential polynomial $(A^t)_{ij} = 0$ when $t=0$ and $i \neq j$.  For each integer $k$, if the interval $(k,k+1)$ contains an exponent $t$ such that $(A^t)_{ij}$ is negative, then the closed interval $[k,k+1]$ must contain at least two zeros of $(A^t)_{ij}$ since $A^k$ and $A^{k+1}$ are nonnegative.  If $w_{ij} = 3$ or $4$, then there is at most one integer $k > 0$ such that the interval $(k,k+1)$ contains an exponent $t$ with $(A^t)_{ij} < 0$, otherwise there would be more than 4 zeros counting multiplicity.  

Let $s$ be any real number.  Note that if every entry of a row of $A^t$ is nonnegative for all $t \in (s,s+1)$, then the row will continue to be nonnegative for all $t > s$.  This is because the rows of $A^{t+1}$ are equal to the rows of $A^t$ multiplied by the nonnegative matrix $A$ on the right.  A similar observation applies to the columns of $A^t$.  Thus if any row of $A^t$ (or column) contains  negative entry for some $t \in \reals$, then that same row (or column) of $A^{t-k}$ contains a negative entry for all integers $k > 0$.  

Suppose that row $i$ of $W$ contains $M$ entries with values greater than $2$ and none with values greater than 4.  If row $i$ of $A^t$ has an entry that is negative for some $t > M$, then it must have negative entries for exponents $t$ in each of the intervals $(1,2)$, $(2,3)$, ..., $(M,M+1)$.  Each of these $M$ intervals must then correspond uniquely to one of the $M$ entries of row $i$ of $W$ that are larger than 2.  Since there are no other entries with critical exponents larger than 1, we conclude that that row $i$ of $A^t$ has no negative entries for any $t > M+1$.  The proof for a column is identical.  
\end{proof}

\begin{proof}[Proof of Theorem \ref{5by5}]
By continuity, it suffices to prove the result in the generic case where $A$ is invertible, $A$ has 5 distinct eigenvalues $\lambda_1, \ldots, \lambda_5$, and the eigenvectors of $A$ have no zero entries.  For any such $A$, the corresponding sign change matrix must be one of the examples listed below (up to permutation similarity).  This list of possible 5-by-5 sign change matrices was generated using MATLAB by taking all possible sign patterns of the eigenvector matrix $U$, with the requirements that one vector be positive (Perron), that the top entry of each column be positive, and that no two rows or columns have the same sign pattern.  

$$
\begin{bmatrix} 
0 & 1 & 2 & 2 & 2 \\ 
1 & 0 & 1 & 3 & 3 \\ 
2 & 1 & 0 & 2 & 4 \\ 
2 & 3 & 2 & 0 & 2 \\ 
2 & 3 & 4 & 2 & 0 \\
\end{bmatrix}, ~
\begin{bmatrix} 
0 & 1 & 2 & 2 & 3 \\ 
1 & 0 & 1 & 3 & 2 \\ 
2 & 1 & 0 & 2 & 3 \\ 
2 & 3 & 2 & 0 & 3 \\ 
3 & 2 & 3 & 3 & 0 \\ 
\end{bmatrix}, ~
\begin{bmatrix} 
0 & 1 & 2 & 3 & 2 \\ 
1 & 0 & 1 & 2 & 3 \\ 
2 & 1 & 0 & 1 & 4 \\ 
3 & 2 & 1 & 0 & 3 \\ 
2 & 3 & 4 & 3 & 0 \\ 
\end{bmatrix}
$$
$$
\begin{bmatrix} 
0 & 1 & 2 & 3 & 3 \\ 
1 & 0 & 1 & 2 & 2 \\ 
2 & 1 & 0 & 1 & 3 \\ 
3 & 2 & 1 & 0 & 2 \\ 
3 & 2 & 3 & 2 & 0 \\ 
\end{bmatrix}, ~ 
\begin{bmatrix} 
0 & 1 & 2 & 2 & 3 \\ 
1 & 0 & 1 & 3 & 4 \\ 
2 & 1 & 0 & 2 & 3 \\ 
2 & 3 & 2 & 0 & 1 \\ 
3 & 4 & 3 & 1 & 0 \\ 
\end{bmatrix},~
\begin{bmatrix} 
0 & 1 & 2 & 3 & 4 \\ 
1 & 0 & 1 & 2 & 3 \\ 
2 & 1 & 0 & 1 & 2 \\ 
3 & 2 & 1 & 0 & 1 \\ 
4 & 3 & 2 & 1 & 0 \\ 
\end{bmatrix}$$ 

$$\begin{bmatrix} 
0 & 1 & 2 & 2 & 2 \\ 
1 & 0 & 1 & 1 & 3 \\ 
2 & 1 & 0 & 2 & 4 \\ 
2 & 1 & 2 & 0 & 2 \\ 
2 & 3 & 4 & 2 & 0 \\ 
\end{bmatrix},~
\begin{bmatrix} 
0 & 1 & 2 & 2 & 3 \\ 
1 & 0 & 1 & 1 & 2 \\ 
2 & 1 & 0 & 2 & 3 \\ 
2 & 1 & 2 & 0 & 1 \\ 
3 & 2 & 3 & 1 & 0 \\ 
\end{bmatrix},~ 
\begin{bmatrix} 
0 & 1 & 2 & 2 & 3 \\ 
1 & 0 & 1 & 1 & 4 \\ 
2 & 1 & 0 & 2 & 3 \\ 
2 & 1 & 2 & 0 & 3 \\ 
3 & 4 & 3 & 3 & 0 \\ 
\end{bmatrix} $$

$$\begin{bmatrix} 
0 & 1 & 2 & 2 & 2 \\ 
1 & 0 & 1 & 3 & 3 \\ 
2 & 1 & 0 & 2 & 2 \\ 
2 & 3 & 2 & 0 & 2 \\ 
2 & 3 & 2 & 2 & 0 \\ 
\end{bmatrix},~
\begin{bmatrix} 
0 & 1 & 2 & 2 & 2 \\ 
1 & 0 & 1 & 1 & 1 \\ 
2 & 1 & 0 & 2 & 2 \\ 
2 & 1 & 2 & 0 & 2 \\ 
2 & 1 & 2 & 2 & 0 \\ 
\end{bmatrix},~ 
\begin{bmatrix} 
0 & 1 & 2 & 1 & 4 \\ 
1 & 0 & 1 & 2 & 3 \\ 
2 & 1 & 0 & 1 & 2 \\ 
1 & 2 & 1 & 0 & 3 \\ 
4 & 3 & 2 & 3 & 0 \\ 
\end{bmatrix} $$

$$\begin{bmatrix} 
0 & 1 & 3 & 2 & 2 \\ 
1 & 0 & 2 & 3 & 3 \\ 
3 & 2 & 0 & 1 & 3 \\ 
2 & 3 & 1 & 0 & 2 \\ 
2 & 3 & 3 & 2 & 0 \\ 
\end{bmatrix},~
\begin{bmatrix} 
0 & 1 & 3 & 2 & 3 \\ 
1 & 0 & 2 & 3 & 4 \\ 
3 & 2 & 0 & 1 & 2 \\ 
2 & 3 & 1 & 0 & 1 \\ 
3 & 4 & 2 & 1 & 0 \\ 
\end{bmatrix},~ 
\begin{bmatrix} 
0 & 1 & 2 & 4 & 3 \\ 
1 & 0 & 1 & 3 & 4 \\ 
2 & 1 & 0 & 2 & 3 \\ 
4 & 3 & 2 & 0 & 1 \\ 
3 & 4 & 3 & 1 & 0 \\ 
\end{bmatrix} $$

$$\begin{bmatrix} 
0 & 1 & 3 & 2 & 3 \\ 
1 & 0 & 2 & 3 & 4 \\ 
3 & 2 & 0 & 3 & 2 \\ 
2 & 3 & 3 & 0 & 1 \\ 
3 & 4 & 2 & 1 & 0 \\ 
\end{bmatrix},~ 
\begin{bmatrix} 
0 & 1 & 3 & 3 & 3 \\ 
1 & 0 & 2 & 4 & 2 \\ 
3 & 2 & 0 & 2 & 2 \\ 
3 & 4 & 2 & 0 & 2 \\ 
3 & 2 & 2 & 2 & 0 \\ 
\end{bmatrix},~ 
\begin{bmatrix} 
0 & 1 & 2 & 4 & 2 \\ 
1 & 0 & 1 & 3 & 3 \\ 
2 & 1 & 0 & 2 & 4 \\ 
4 & 3 & 2 & 0 & 2 \\ 
2 & 3 & 4 & 2 & 0 \\ 
\end{bmatrix}$$

$$\begin{bmatrix} 
0 & 1 & 2 & 3 & 3 \\ 
1 & 0 & 1 & 4 & 2 \\ 
2 & 1 & 0 & 3 & 3 \\ 
3 & 4 & 3 & 0 & 2 \\ 
3 & 2 & 3 & 2 & 0 \\ 
\end{bmatrix} ,~
\begin{bmatrix} 
0 & 1 & 3 & 2 & 3 \\ 
1 & 0 & 2 & 3 & 2 \\ 
3 & 2 & 0 & 3 & 2 \\ 
2 & 3 & 3 & 0 & 3 \\ 
3 & 2 & 2 & 3 & 0 \\ 
\end{bmatrix} ,~
\begin{bmatrix} 
0 & 1 & 3 & 3 & 3 \\ 
1 & 0 & 2 & 2 & 2 \\ 
3 & 2 & 0 & 2 & 2 \\ 
3 & 2 & 2 & 0 & 2 \\ 
3 & 2 & 2 & 2 & 0 \\ 
\end{bmatrix} $$

For this list of possible sign change matrices, it is straightforward to verify that every entry is contained in either a row or a column with fewer than three entries greater than 2 and none greater than 4.  Therefore we can apply Lemmas \ref{Wbound12} and \ref{Wbound34} to see that the critical exponent for any 5-by-5 invertible $\DN$ matrix $A$ is at most 3.  
\end{proof}

\section{Additional Observations}

\begin{theorem}
Let $A$ be an $n$-by-$n$ doubly nonnegative matrix.  Then there is an $\epsilon>0$ such that $(\epsilon A+I)^t$ is doubly nonnegative for all $t \geq n-2$.   
\end{theorem}
\begin{proof}
Suppose that $A$ is irreducible.  By Lemma \ref{primitive}, $A^{n-1}$ has all positive entries, as does $A^n$.  Therefore we may choose $\epsilon > 0$ small enough that $\epsilon A^n \leq A^{n-1}$ entry-wise.  Then $(\epsilon A)^k \leq (\epsilon A)^{k-1}$ for all $k \geq n$.  
We now express $(\epsilon A+I)^t$ using the binomial series:
$$(\epsilon A+I)^t = \sum_{k = 0}^{\infty} c_k (\epsilon A)^k$$
where 
$$c_k =\left( \frac{t (t-1) (t-2) \cdots (t-k +1)}{k!} \right).$$
If $t > n-2$, then $c_k \geq 0$ for all $k \leq n-1$.  Furthermore, for all $k > t$, $|c_k| < |c_{k-1}|$.  In particular, if $c_k$ is negative, then $|c_{k}| < |c_{k-1}|$.  Together with the fact that $(\epsilon A)^k$ is entry-wise decreasing for $k \geq n-1$, we conclude that $(\epsilon A+I)^t$ has all positive entries for all $t > n-2$.  

Since $A$ is symmetric, if $A$ is reducible, then there is a permutation matrix $P$ such that $PAP^T$ is the direct sum of irreducible matrices or possibly 1-by-1 blocks containing zero.  Each of these blocks will be nonnegative when raised to any power greater than or equal to $n-2$, so the same applies to $A$.  
\end{proof}

\begin{theorem}
Suppose $A$ is a doubly nonnegative matrix such that the largest eigenvalue $\lambda_1$ has multiplicity one and corresponding eigenvector $x_1$.  Then for any constant $r > 0$, the matrix $B = A + r x_1 x_1^T$ satisfies $B^t \geq A^t$ entry-wise for all $t \geq 0$.  That is, the entries of $A^t$ are monotone with respect to the largest eigenvalue.   
\end{theorem}
\begin{proof}
By the Perron-Frobenius theorem, $x_1$ has all nonnegative entries.  The result of the theorem follows directly.  
\end{proof}

Although we have conjectured that the critical exponent for $n$-by-$n$ doubly nonnegative matrices is $n-2$, an even stronger result may be true.  The critical exponent of a doubly nonnegative matrix might only depend on the number of distinct eigenvalues.  The following theorem shows that for doubly nonnegative matrices with only three distinct eigenvalues, the critical exponent is at most $1$.  This suggests the following conjecture: for an $n$-by-$n$ doubly nonnegative matrix with $k$ distinct eigenvalues, the critical exponent may be at most $k-2$.   

\begin{theorem} If $A$ is a doubly nonnegative matrix with at most three distinct eigenvalues, then $A^t$ is doubly nonnegative for all $t\geq 1$.  
\end{theorem}
\begin{proof}
By continuity, it suffices to assume that $A$ is invertible and therefore that $\lim_{t \rightarrow 0} A^t = I_n$.  Since $A$ has at most three distinct eigenvalues, $\lambda_1, \lambda_2,$ and $\lambda_3$, note that the exponential polynomials for each entry $(A^t)_{ij}$ can have at most 2 sign changes.  Thus, the exponential polynomial for each entry of $A^t$ can have at most 2 zeros.  For off-diagonal entries, one of the zeros is $t=0$.  Since there will be only one remaining zero counting multiplicity, and $A^t$ is nonnegative for $t=1$, the off diagonal entries must be positive for all $t>1$.  Since $A$ is positive definite, the diagonal entries of $A^t$ are positive for all $t$.  
\end{proof}

\bibliography{DNrefs}

\begin{thebibliography}{1}

\bibitem{Chen}
Shencan Chen.
\newblock Proof of a conjecture concerning the {H}adamard powers of inverse
  {$M$}-matrices.
\newblock {\em Linear Algebra Appl.}, 422(2-3):477--481, 2007.

\bibitem{FH}
Carl~H. FitzGerald and Roger~A. Horn.
\newblock On fractional {H}adamard powers of positive definite matrices.
\newblock {\em J. Math. Anal. Appl.}, 61(3):633--642, 1977.

\bibitem{HornJohnson1}
Roger~A. Horn and Charles~R. Johnson.
\newblock {\em Matrix analysis}.
\newblock Cambridge University Press, Cambridge, 1985.

\bibitem{HornJohnson2}
Roger~A. Horn and Charles~R. Johnson.
\newblock {\em Topics in matrix analysis}.
\newblock Cambridge University Press, Cambridge, 1994.
\newblock Corrected reprint of the 1991 original.

\bibitem{PoSz}
George P{\'o}lya and Gabor Szeg{\H{o}}.
\newblock {\em Problems and theorems in analysis. {II}}.
\newblock Classics in Mathematics. Springer-Verlag, Berlin, 1998.
\newblock Theory of functions, zeros, polynomials, determinants, number theory,
  geometry, Translated from the German by C. E. Billigheimer, Reprint of the
  1976 English translation.

\end{thebibliography}
\bibliographystyle{plain}
\end{document}